\documentclass[10pt,reqno]{amsart}

\usepackage{amsmath,amsfonts,amsbsy,amsgen,amscd,mathrsfs,amssymb,amsthm}
\usepackage{enumerate,mathtools}

\usepackage[usenames,dvipsnames]{xcolor}
\usepackage[colorlinks=true,citecolor=blue,linkcolor=blue]{hyperref}
\usepackage{bbm}
\usepackage{tikz}
\usetikzlibrary{shadings}

\newtheorem{thm}{Theorem}[section]
\newtheorem{theorem}{Theorem}
\newtheorem{corollary}{Corollary}
\newtheorem{lem}[thm]{Lemma}

\theoremstyle{definition}

\newtheorem{rem}[thm]{Remark}

\numberwithin{equation}{section} 
\numberwithin{figure}{section}
\numberwithin{table}{section}

\begin{document}

\title{Sharpening the gap between $L^{1}$ and $L^{2}$ norms}

\author{Paata Ivanisvili}
\address{(P.I.) Department of Mathematics, University of California, 
Irvine, CA 92617, USA}
\email{pivanisv@uci.edu}

\author{Yonathan Stone}
\address{(Y.S.) Department of Mathematics, University of California, 
Irvine, CA 92617, USA}
\email{ystone@uci.edu}

\begin{abstract}
We refine the  classical Cauchy--Schwartz inequality $\|X\|_{1} \leq \|X\|_{2}$ by demonstrating that for any 
 $p$ and $q$ with $q>p>2$,  there exists a constant $C=C(p,q)$ such that 
 \begin{align*}
  \|X\|_1 \leq 1 - C \cdot \frac{\Big{(}\|X\|_p^p - 1\Big{)}^{\frac{q-2}{q-p}}}{\Big{(}\|X\|_q^q - 1\Big{)}^{\frac{p-2}{q-p}}}
 \end{align*}
 holds true for all Borel measurable random variables $X$ with $\|X\|_{2}=1$ and  $\|X\|_{p}<\infty$. We illustrate two applications of this result: one for  biased Rademacher sums and another for exponential sums. 
\end{abstract}


\maketitle

\thispagestyle{empty}

\section{Introduction}
\subsection{Separating first and second moments}

Let $\mathcal{F}$ be a family of  Borel measurable random variables. For any $p>0$, define $\| X\|_{p} := (\mathbb{E} |X|^{p})^{1/p}$. Assuming  $\|X\|_{2}=1$ for any $X \in \mathcal{F}$, then the Cauchy--Schwarz inequality tells us that $\|X\|_{1}\leq 1$. In this paper, our goal is to understand under what conditions on $X$ there exists $\varepsilon >0$ such that 
\begin{align}\label{uketesi}
\|X\|_{1} \leq 1-\varepsilon
\end{align}
 holds for all $X \in \mathcal{F}$. 

The conditions on $X$ should be such that one can easily verify them. There are two typical examples we keep in mind: $\mathcal{F} = \{  \sum_{j=1}^{n}e^{2 \pi j^{2} i \theta }, \,  n\geq 2\}$ with $\theta \sim  \mathrm{unif}([0,1])$, and $\mathcal{F} =\{ \sum_{j=1}^{n} a_{j}\xi_{j}, \, n\geq 1, \, a_{j} \in \mathbb{R}\}$, where $\xi_{j}$ are i.i.d. random variables. Notice that if $\|X\|_{1}\leq 1-\varepsilon$, then it follows from H\"older's inequality that $1 \leq \|X\|_{1}^{\theta} \|X\|_{p}^{1-\theta}$ for any $p>2,$ where $\theta = \frac{p-2}{2(p-1)}$. Hence, we arrive at the necessary condition $\|X\|_{p} \geq 1+C(p)\varepsilon$  for some positive constant $C(p)>0$, which in practice can be verified for specific values of $p$. However, random variables $X$ taking two values show that this condition is not sufficient for (\ref{uketesi}). 

The main result in this paper shows that if the necessary condition holds for some $p>2$ and we also have good control on the growth of $\|X\|_{q}$ for some $q>p$, then we get (\ref{uketesi}).

\begin{theorem}\label{main2}
    Let $q > p > 2$ be finite. Then there exists a constant $C=C(p,q)>0$ such that 
    \begin{equation*}
       \|X\|_1 \leq 1- C \frac{\Big{(}\|X\|_p^p - 1\Big{)}^{\frac{q-2}{q-p}}}{\Big{(}\|X\|_q^q - 1\Big{)}^{\frac{p-2}{q-p}}}
    \end{equation*}
    holds for any Borel measurable random variable $X$ with $\|X\|_{2}=1$ and  $\|X\|_{p}<\infty$.  
\end{theorem}
\begin{rem}
We will see from the proof that
\begin{align*}
&C(p,q)=\\
&\inf_{a,c \in (0,1)} \frac{(c^{q-2}(a^{q}-1)+c^{q}(a^{2}-a^{q})+1-a^{2})^{\frac{p-2}{q-p}}}{(c^{p-2}(a^{p}-1)+c^{p}(a^{2}-a^{p})+1-a^{2})^{\frac{q-2}{q-p}}}\cdot (1-c)(1-a)(1-ac),
\end{align*}
and this result is sharp, i.e., the sharpness of the constant $C(p,q)$ can be verified using random variables taking on two values. Our proof will show 
$$
C(p,q) \geq \frac{(\min\{1,q-2\})^{\frac{p-2}{q-p}}}{p^{\frac{2(q-2)}{q-p}}} > 0
$$
however, sometimes the constant $C(p,q)$ can be explicitly computed for given powers $p$ and $q$. For example, we will see that $C(4,6) = 1/3$. 
\end{rem}

\begin{rem}
One could use arguments from the classical moment problem, specifically the positive semidefiniteness of two matrices whose entries are the moments of \(X\) chosen in a particular way (see page 781 in \cite{BP}), to prove Theorem~\ref{main2} when \(p\) and \(q\) are integers. This approach seems feasible for small integers \(p\) and \(q\), such as \(p = 4\) and \(q = 6\). However, the computational complexity may increase with larger values of \(p\) and \(q\). Our method is different and does not require \(p\) and \(q\) to be integers.
\end{rem}

\subsection{An application to exponential sums}

In \cite{Bu1} Bourgain asked the following question: does there exist $\varepsilon>0$ such that for any finite subset $S \subset \mathbb{\mathbb{Z}}$ of distinct integers with $|S|>1$ we have 
\begin{align}\label{bou2}
\int_{0}^{1} \left|\frac{1}{\sqrt{|S|}} \sum_{j \in S} e^{2\pi i j \theta}\right|d\theta \leq 1-\varepsilon. 
\end{align}
It was proved in \cite{CA1} that if such $\varepsilon >0$ exists then it must be at most $1-\frac{\sqrt{\pi}}{2}$. Denote 
$$
X_{S} = \frac{1}{\sqrt{|S|}} \sum_{j\in S} e^{2\pi i j \theta},
$$
 where $\theta \sim \mathrm{unif}[0,1]$. The best result in (\ref{bou2})  is due to Bourgain \cite{Bu1}, who showed existence of a constant $c>0$ such that  
\begin{align*}
\| X_{S}\|_{1} \leq  1-c\frac{\log(|S|)}{|S|}.
\end{align*}

Bourgain's question is already nontrivial  for squares of integers, i.e.,  
$$
Q := \{ j^{2}\,  |\,  j=1, \ldots, m\},
$$
  and to the best of our knowledges the inequality (\ref{bou2}) is open in this case. In this case Bourgain~\cite{Bu3} showed that for any $q >4$ we have 
\begin{align}
\|X_{Q}\|^{q}_{q} &\leq  C(q) |Q|^{\frac{q}{2}-2}  \quad \text{for some}  \quad C(q)<\infty, \label{jean1}\\
\| X_{Q} \|^{4}_{4} &\geq C \log(|Q|) \quad  \text{with some positive $C>0$}. \label{jean2}
\end{align}  
Combining the estimates (\ref{jean1}) and (\ref{jean2}) together with Theorem~\ref{main2} applied with $p=4$ we obtain
\begin{align*}
\| X_{Q}\|_{1} \leq 1-C'(q) \frac{\log^{\frac{q-2}{q-4}}(|Q|)}{|Q|}
\end{align*} 
holds with some constant $C'(q)>0$. Given any $N>1$ we can choose  $q>4$  so that $\frac{q-2}{q-4}=N$, thus we obtain
\begin{corollary}
For any  $N>0$ there exists a constant $c(N)>0$ such that 
\begin{align*}
\| X_{Q}\|_{1} \leq  1-c(N)\frac{\log^{N}(|Q|)}{|Q|}. 
\end{align*}
\end{corollary}

\subsection{An application to $L^{1}$ Poincar\'e inequality on the hypercube}

Fix $ p \in (0,1)$, and let $\xi_{1}, \ldots, \xi_{n}$ be i.i.d. Bernoulli random variables such that 
\begin{align}\label{biased1}
\xi_{1} = 
\begin{cases}
&\sqrt{\frac{1-p}{p}} \quad \text{with probability} \quad p,\\
&-\sqrt{\frac{p}{1-p}}\quad \text{with probability} \quad 1-p.
\end{cases}
\end{align}
Clearly $\mathbb{E} \xi_{1}=0$, and $\mathbb{E} |\xi_{1}|^{2}=1$. The following theorem was proved in \cite{pi1}. 
\begin{theorem}
There exists $\varepsilon>0$ such that 
\begin{align*}
\sup_{n \geq 1} \sup_{a_{1}^{2}+\ldots+a_{n}^{2}=1} \mathbb{E} \left|\sum_{j=1}^{n}a_{j} \xi_{j}\right|<1
\end{align*}
holds for any $p \in (3/4 - \varepsilon, 3/4+\varepsilon)$. 
\end{theorem}

Notice that the conclusion of the theorem does not hold  for $p=1/2$. The theorem was one of the technical steps in proving existence of a small constant $\delta>0$ such that 
\begin{align}\label{poinc}
\mathbb{E} |f(x)-\mathbb{E} f(x)|\leq (\frac{\pi}{2}-\delta) \mathbb{E} |\nabla f|(x)
\end{align}
holds for any $f :\{-1,1\}^{n} \mapsto \mathbb{R}$, all $n\geq 1$, where $x=(x_{1}, \ldots, x_{n}) \sim \{-1,1\}^{n}$, 
\begin{align*}
|\nabla f|(x) = \left( \sum_{j=1}^{n} |D_{j}f(x)|^{2}\right)^{1/2},
\end{align*}
and 
\begin{align*}
D_{j}f(x) := \frac{f(x)-f(S_{j}(x))}{2}, \quad \text{here} \quad S_{j}(x) = (x_{1}, \ldots, x_{j-1}, -x_{j}, x_{j+1}, \ldots, x_{n}). 
\end{align*}

The estimate (\ref{poinc}) was particularly surprising due to the variety of proofs leading to the same bound $\mathbb{E} |f(x)-\mathbb{E} f(x)| \leq \frac{\pi}{2} \mathbb{E} |\nabla f|(x)$. Notably, one of these proofs employed a non-commutative approach by L.~Ben-Efraim and F.~Lust-Piquard \cite{Ben}. The consistent appearance of the constant $\pi/2$ across these proofs suggested that the Cheeger constant in the $L^{1}$ Poincar\'e inequality might be $\pi/2$. However, this conjecture was refuted in \cite{pi1} by establishing (\ref{poinc}).

Recently R.~van Handel \cite{rvan} obtained more quantitative bound on $\delta$ by showing  
\begin{align*}
\|X\|_{1} \leq \sqrt{1-\frac{(\|X\|_{4}^{4}-1)^{2}}{32 \|X\|_{6}^{6}}}
\end{align*}
hold for all random variables $X$ with $\|X\|_{2}=1$, and $\|X\|_{4}<\infty$, which allowed him to obtain
\begin{align}\label{ramon1}
\sup_{n \geq 1} \sup_{a_{1}^{2}+\ldots+a_{n}^{2}=1} \mathbb{E} \left|\sum_{j=1}^{n}a_{j} \xi_{j}\right|\leq \sqrt{1- \frac{\min\left\{4p^{3}(1-p)^{3},(2p-1)^{4}p(1-p)\right\}}{480}}.
\end{align}
The inequality (\ref{ramon1}) combined together with techniques from \cite{pi2} gives $\delta \approx 0.000006...$.

As an immediate application of Theorem~\ref{main2} we will show the following 
\begin{corollary}\label{texnika1}
We have 
\begin{align}\label{dax1}
\sup_{n \geq 1} \sup_{a_{1}^{2}+\ldots+a_{n}^{2}=1} \mathbb{E} \left|\sum_{j=1}^{n}a_{j} \xi_{j}\right|\leq 1-\frac{\min\left\{4p^{3}(1-p)^{3}, (2p-1)^{4}p(1-p)\right\}}{45-3(p(1-p))^{3}}
\end{align}
holds for all $p \in [0,1]$. 
\end{corollary}
The corollary combined with arguments from \cite{pi2} gives the following 
\begin{corollary}\label{texnika2}
For any $f :\{-1,1\}^{n} \mapsto \mathbb{R}$ we have that
\begin{align}\label{puank1}
\mathbb{E} |f(x) - \mathbb{E} f(x)| \leq \left(\frac{\pi}{2}-\delta\right) \mathbb{E} |\nabla f|(x)
\end{align}
holds with $\delta \approx 0.00013...$
\end{corollary}

\begin{rem}
It seems to us that the best  $\delta$ one can choose in (\ref{puank1}) is  $\frac{\pi}{2} - \sqrt{\frac{\pi}{2}} \approx 0.31748...$, at least this is the case for an analogous question in Gauss space, however, to prove such a statement one needs to come up with techniques different from Section~\ref{datvla}. Indeed, let $X = \sum_{j=1}^{n} a_{j}\xi_{j}$ be a biased Rademacher sum. Let $R=R(p) \in (0,1)$ be the best constant,  for each $p \in (0,1)$,  such that  $\| X\|_{1}/\|X\|_{2} \leq R(p)$. In \cite{yonathan} it was proved that 
\begin{align}\label{stone}
\sup_{n\geq 1, \; a^{2}_{j}=a_{j}} \frac{\| X\|_{1}}{\|X\|_{2}} = 2 \sqrt{N(p^{*})p^{*}} (1-p^{*})^{N(p^{*})-1},
\end{align}
where $N(p) = \lfloor \frac{1}{1-(1-p)^{2}} \rfloor$, and $p^{*}=\min(p, 1-p)$. Therefore, the best $\delta$ one may hope to get in (\ref{puank1}) using the argument in Section~\ref{datvla}, assuming $R(p)$ equals to the right  hand side of (\ref{stone}), is $\frac{\pi}{2} - \int_{1/2}^{1}2 \sqrt{N((1-p))(1-p)} p^{N(1-p)-1} \frac{dp}{\sqrt{p(1-p)}} \approx 0.149...$ which is still far from $\pi/2 - \sqrt{\pi/2} \approx 0.31748...$. 
\end{rem}

\section{Proofs}

\subsection{Proof of Theorem \ref{main2}}

By rescaling we can rewrite the claimed inequality in the theorem as follows
\begin{align}\label{ineq2}
 \|X\|_2 - \|X\|_1 \geq \frac{C(p,q)}{\|X\|_2} \cdot \frac{\Big{(}\|X\|_p^p - \|X\|_2^p\Big{)}^{\frac{q-2}{q-p}}}{\Big{(}\|X\|_q^q - \|X\|_2^q\Big{)}^{\frac{p-2}{q-p}}}
 \end{align}
 which holds for all Borel random variables $X$ with $\|X\|_{p}<\infty$.  If $\|X\|_{q} = \infty$ there is nothing to prove, so in what follows we consider random variables $X$ such that  $\| X\|_{q}<\infty$.  Without loss of generality, we will first assume that $\|X\|_\infty < \infty$ since by the Dominated Convergence Theorem, we can take $X_N = \min\{|X|, N\}$ and note that $X_N \leq |X|$ and 
    \[X_N \rightarrow |X| \text{\quad a.e.}\]
    Furthermore, by homogeneity of (\ref{ineq2}), assume for now that $0\leq X \leq 1$.  Note that since we take absolute values in all $L^p$ norms that it suffices to prove (\ref{ineq2}) when $X \geq 0$. As for the assumption that $X \leq 1$, we will briefly need this assumption, but note that we will soon replace it with a more useful stipulation.  We now consider the following optimization problem:
    \[\Phi(x,y,z) = \sup\limits_{0 \leq X \leq 1} \{ \mathbb{E}X,\quad \mathbb{E}X^2 = x, \quad \mathbb{E}X^p = y, \quad \mathbb{E}X^q = z\}.\]
    We will note that fortunately for us the space curve $\gamma(t) = (t^2, t^p, t^q, -t),\quad t\in [0,1]$ has totally positive torsion.  We recall that this means that the leading principal minors of the $4\times 4$ matrix $(\gamma^{(1)}(t), \gamma^{(2)}(t), \gamma^{(3)}(t), \gamma^{(4)}(t))$ are positive for all $t \in (0,1)$.  For the sake of completeness, we will check  total positivity here.  First of all, the matrix is given by
    \[\begin{bmatrix*}
        \gamma^{(1)}(t) \\ \gamma^{(2)}(t) \\ \gamma^{(3)}(t) \\ \gamma^{(4)}(t)
    \end{bmatrix*} = \begin{bmatrix*}[l]
        2t & pt^{p-1} & qt^{q-1} & -1 \\ 
        2 & p(p-1)t^{p-2} & q(q-1)t^{q-2} & 0 \\
        0 & p(p-1)(p-2)t^{p-3} & q(q-1)(q-2)t^{q-3} & 0 \\
        0 & p(p-1)(p-2)(p-3)t^{p-4} & q(q-1)(q-2)(q-3)t^{q-4} & 0
    \end{bmatrix*}\]
    and that the $i$-th leading principal minors $A_{ii}$ are each respectively given by
    \begin{align*}
        A_{11} &= 2t, \\
        A_{22} &= t^{p-1}2p(p-1),\\
        A_{33} &= 2pq(p-2)(q-2)(q-p)t^{p+q-4}, \\
        A_{44} &= 2p(p-1)(p-2)q(q-1)(q-2)(q-p)t^{p+q-7},
    \end{align*}
    where positivity follows from the fact that $t > 0$ and that $2 < p < q$.  We can thus invoke Theorem 2.3, case 1 in \cite{pi4}, which tells us that the supremum in $\Phi(x,y,z)$ is achieved by a random variable $X$ taking on at most two values.  To that end, we will assume in the sequel that
    \[\mathbb{P}\{X = a\} = 1 - \mathbb{P}\{X=b\} = r.\]
    As we no longer need it, we will drop the assumption $X \leq 1$ at this point and by homogeneity instead replace it with the assumption that $\|X\|_2 = 1$.  Note that since
    \[\|X\|_2^2 = a^2r + b^2(1-r),\]
    our new normalization assumption implies that
    \[r = \frac{b^2 - 1}{b^2 - a^2}.\]
    Taking this into account, we obtain that
    \[\|X\|_1 = \frac{1 + ab}{a + b},\]
meaning the desired inequality can now be expressed as
\[\frac{(1-a)(b-1)}{a + b} \geq C(p,q) \frac{\left(\tfrac{a^p(b^2-1) + b^p(1-a^2)}{b^2 - a^2}-1\right)^\theta}{\left(\tfrac{a^q(b^2-1) + b^q(1-a^2)}{b^2 - a^2}-1\right)^{\theta - 1}},\]
where $\theta = \frac{q-2}{q-p}$, and we may furthermore assume without loss of generality that $0 < a < 1 < b < \infty$, otherwise arbitrary.  We next note that we can further simplify this as
\begin{align}\label{simp} 
\left(\tfrac{a^q(b^2-1) + b^q(1-a^2)+a^{2}-b^{2}}{(1-a)(b-1)(b-a)}\right)^{\theta - 1} \geq C(p,q)\left(\tfrac{a^p(b^2-1) + b^p(1-a^2)+a^{2}-b^{2}}{(1-a)(b-1)(b-a)}\right)^{\theta}. 
\end{align}
After a change of a variables $b=1/c$, where $c \in (0,1)$, the inequality (\ref{simp}) rewrites and simplifies  as 
\begin{align*}
&\left(\frac{c^{q-2}(a^q-1) + c^q(a^2-a^q) + 1-a^2 }{(1-c)(1-a)(1-ac)}\right)^{\theta - 1} \geq \\
&C(p,q)\left(\frac{c^{p-2}(a^p-1) + c^p(a^2-a^p) + 1-a^2 }{(1-c)(1-a)(1-ac)}\right)^{\theta}. 
\end{align*}
Next,  by defining
\[B(a,c,q) := \frac{c^{q-2}(a^q-1) + c^q(a^2-a^q) + 1-a^2 }{(1-c)(1-a)(1-ac)},\]
we can take 
\begin{align}\label{znachc}
C(p,q) = \inf_{a,c \in (0,1)} \frac{B(a,c,q)^{\theta-1}}{B(a,c,p)^{\theta}}.
\end{align}
To show $C(p,q)>0$ we claim that  we can find constants $c_1(p), c_2(p)> 0 $ such that 
\[c_1(p) \leq B(a,c,p) \leq c_2(p),\]
Then the equality (\ref{znachc}) will imply
\[C(p,q) > \frac{c_1(q)^{\theta-1}}{c_2(p)^\theta}.\]
We will thus move on to finding upper and lower bounds on the function $B(a,c,p)$, where $2 < p < \infty$ and $(a,c) \in (0,1)\times (0,1).$
 \begin{lem}\label{texnika}
 We have $B(a,c,p) \geq \min\{1, p-2\}$ for all $a,c \in (0,1)$ and all $p> 2$. 
 \end{lem}
 \begin{proof}
 Let $L:=\frac{p-2}{p}$, and for each $c \in (0,1)$ consider the family of functions $f(a) = f(a; c)$ defined as 
 \begin{align*}
 f(a) := c^{p-2}(a^p-1) + c^p(a^2-a^p) + 1-a^2 -L(1-c)(1-a)(1-ac).
 \end{align*}
 Clearly the lemma is the same as $f(a) \geq 0$ for all $a, c \in (0,1)$. To verify $f(a) \geq 0$ we will argue as follows: we will show that 
\begin{itemize}
\item[(1)] $f(0)\geq 0$. 
\item[(2)] $f(1)=0$. 
\item[(3)] $f'(1)\leq 0$. 
\item[(4)] $f''(0)<0$. 
\item[(5)] $f''$ changes sign at most once from $-$ to $+$.  
\end{itemize}
It then follows that $f(a) \geq 0$ on $[0,1]$. Indeed, if $f''\leq 0$ on $[0,1]$ then we are done because of (1) and (2). 
If $f''$ changes sign from $-$ to $+$ at a point $k \in [0,1)$, then (2) and (3) imply that $f \geq 0$ on $[k,1]$. In particular $f(k)\geq 0$. 
Thus $f\geq 0$ on $[0,k]$ because $f'' \leq 0$ on $[0,k]$, and this finishes the proof of Lemma~\ref{texnika}.

To verify (1), we have 
\begin{align*}
f(0)=-c^{p-2}+1-L(1-c).
\end{align*} 
Since $\varphi(c) = -c^{p-2}+1-L(1-c)$ has the properties $\varphi(0)=1-L\geq0$, $\varphi(1)=0$, and $\varphi'(c) = -(p-2)c^{p-3}+L$ 
changes sign at most once from $+$ to $-$ if  $p>3$,  and $\varphi'\leq 0$ if $p \in (2,3]$ it follows that $\varphi(c) \geq 0$ on $[0,1)$. 

The verification of (2) is trivial, so we move to verifying (3). We have 
\begin{align*}
f'(a) &= pc^{p-2}a^{p-1}+2c^{p}a-pc^{p}a^{p-1}-2a+L(1-c)(1-ac)+Lc(1-c)(1-a). \\
f'(1)&=pc^{p-2}+c^{p}(2-p)-2+L(1-c)^{2}.
\end{align*}
Let $\psi(c) = pc^{p-2}+c^{p}(2-p)-2+L(1-c)^{2}$. We have $\psi(1)=\psi'(1)=0$.  If $p>3$ then we see the coefficients of  the pseudo-polynomial 
\begin{align*}
\psi(c) = c^{p}(2-p)+pc^{p-2}+Lc^{2}-2Lc+L-2
\end{align*}
when arranged according to decreasing order of powers $c$ have the signs $-++--$. Here, we do not know if $p-2\geq 2$ or $p-2<2$ but nevertheless there will be always two sign changes in the coefficients. Therefore by Descartes rule of signs for pseudo-polynomials (see example \#77 on page 46 in \cite{descart}) we obtain that $\psi(c)$ has at most two roots, therefore, $\psi(c)$ does not have roots on $(0,1)$. Since $\psi(0)=L-2<0$ we get $\psi(c)\leq 0$ on $[0,1)$. If $2<p\leq 3$ we have 
\begin{align*}
\psi'(c) &= p(2-p)c^{p-1}+p(p-2)c^{p-3}+2Lc-2L \\
& =p(p-2)c^{p-3}(1-c^{2})-2L(1-c)\\
& = (1-c)(p(p-2)c^{p-3}(1+c)-2L)\\
&> (1-c)(p(p-2)-2L) \geq 0,
\end{align*}
which implies that $\psi(c) \leq 0$ on $[0,1)$.

Next we verify (4). We have 
\begin{align}
f''(a) &= p(p-1)c^{p-2}a^{p-2}+2c^{p}-p(p-1)c^{p}a^{p-2}-2-2Lc(1-c) \nonumber \\
&=a^{p-2}p(p-1)c^{p-2}(1-c^{2})-2(1-c^{p})-2Lc(1-c).\label{prof}
\end{align}
Therefore $f''(0)=-2(1-c^{p})-2Lc(1-c)<0$ since $c \in (0,1)$. 

Finally we verify (5). If follows from (\ref{prof}) that $ a \mapsto f''(a)$ is increasing, and $f''(0)<0$. Therefore $f''$ can change sign at most once from $-$ to $+$ on $(0,1]$. 
 \end{proof}

\begin{lem}
We have $B(a,c,p) \leq p^{2}$ for all $a,c \in (0,1)$ and all $p\geq 2$. 
\end{lem}
\begin{proof}
 Let $M=p^{2}$, and consider 
 \begin{align*}
 g(a) = c^{p-2}(a^p-1) + c^p(a^2-a^p) + 1-a^2 -M(1-c)(1-a)(1-ac).
 \end{align*}
 It suffices to show that $g$ is concave, $g(1)=0$ and $g'(1)\geq 0$. The claim $g(1)=0$ is trivial. To verify $g'(1)\geq 0$ we have 
 \begin{align*}
 g'(1) = c^{p}(2-p)+pc^{p-2}+Mc^{2}-2Mc-2+M.
 \end{align*}
 Let $h(c) = c^{p}(2-p)+pc^{p-2}+Mc^{2}-2Mc-2+M.$ We have $h(0)=M-2\geq 0$, $h(1)=h'(1)=0$. Also $h''(1)=2(M-(p^{2}-2p))>0$. Since the coefficient in front of $c^{p}$ is negative then  there must exist $k>1$ such that $h(k)=0$. The coefficients of $h(c)$ have at most 3 sign changes, therefore by Descartes rule of signs for pseudo-polynomials $h$ cannot have roots in $(0,1)$, so $h(c)\geq 0,$ and hence $g'(1)\geq 0$ on $[0,1]$.

 To verify concavity of $g$ we have 
 \begin{align*}
g''(a)=a^{p-2}p(p-1)c^{p-2}(1-c^{2})-2(1-c^{p})-2Mc(1-c).
 \end{align*}
 Since $a  \mapsto g''(a)$ is increasing in $a$, it suffices to show that $g''(1)\leq 0$. We have 
 \begin{align*}
 g''(1) = -c^{p}(p+1)(p-2)+c^{p-2}p(p-1)+2Mc^{2}-2Mc-2.
 \end{align*}
 If $p>3$ then the pseudo-polynomial  $v(c)=-c^{p}(p+1)(p-2)+c^{p-2}p(p-1)+2Mc^{2}-2Lc-2$ has at most two sign changes in its coefficients, hence at most two roots (counting with multiplicities) on $[0, \infty)$. On the other hand $v(1)=v'(1)=0$, and $v''(1)=-4p(p-2)^{2}<0$, therefore $v(c) \leq 0$ on $[0,1]$. Assume $2<p\leq 3$. Consider $w(c) = v(c)/c$. We have $w(1)=0$. Let us show that $w'(c) \geq 0$ on $(0,1)$. Indeed 
 \begin{align*}
 w'(c)&=p(p-1)(p-3)c^{p-4}-(p^{2}-1)(p-2)c^{p-2}+2M+\frac{2}{c^{2}} \\
 &\geq p(p-1)(p-3)c^{p-4}-(p^{2}-1)(p-2)+2M+\frac{2}{c^{2}}:=N(c).
 \end{align*}
The expression  $N(c)$ satisfies $N(0)=+\infty$ and $N(1)=-2p(p-2)+2M>0$. Also Notice that $N'(c)=(p-4)(p-3)(p-1)pc^{p-5}-\frac{4}{c^{3}}=0$ if and only if 
\begin{align*}
c =c_{0}(p)= \left(\frac{4}{(4-p)(3-p)p(p-1)}\right)^{1/(p-2)}.
\end{align*}
On the other hand we have 
\begin{align*}
c_{0}^{2}N(c_{0}) &= \frac{4(p-3)p(p-1)}{(4-p)(3-p)p(p-1)}+c^{2}(2M-(p^{2}-1)(p-2))+2\\
&=-\frac{2(p-2)}{4-p}+c_{0}^{2}(2M-(p^{2}-1)(p-2)).
\end{align*}
Since $M=p^{2} \geq 2(p-1)(p-2)$ for $p \in [2,3]$, we have 
\begin{align*}
&-\frac{2(p-2)}{4-p}+c_{0}^{2}(2M-(p^{2}-1)(p-2)) \geq \\
&-\frac{2(p-2)}{4-p}+c_{0}^{2}(4(p-1)(p-2)-(p^{2}-1)(p-2))=\\
&-\frac{2(p-2)}{4-p} + \left(\frac{4}{(4-p)(3-p)p(p-1)}\right)^{2/(p-2)}(p-1)(p-2)(3-p)=\\
&\frac{2t}{2-t}\left[ -1+\left(\frac{2^{4-t}}{(2-t)^{2-t}(1-t)^{2-t}(2+t)^{2}(1+t)^{2-t}}\right)^{1/t} \right],
\end{align*}
where $p=2+t$ for $t \in [0,1]$. Consider 
\begin{align*}
f(t) = \log\left(\frac{2^{4-t}}{(2-t)^{2-t}(1-t)^{2-t}(2+t)^{2}(1+t)^{2-t}}\right).
\end{align*}
We have $f(0)=f'(0)=0$. To show $f(t)\geq 0$ on $[0,1]$ it suffices to verify that $f''(t)\geq 0$ on $[0,1]$. A direct computation shows  
\begin{align*}
f''(t) = \frac{-3t^{6}-14t^{5}+8t^{4}+52t^{3}-t^{2}-38t+32}{(2-t)(1-t)^{2}(t+2)^{2}(1+t)^{2}}.
\end{align*}
Notice that 
\begin{align*}
&-3t^{6}-14t^{5}+8t^{4}+52t^{3}-t^{2}-38t+32\geq \\
&-3t^{4}-14t^{4}+8t^{4}+52t^{3}-t^{2}-38t+32=\\
&-9t^{4}+52t^{3}-t^{2}-38t+32\geq 43t^{3}-39t+32.
\end{align*}
On the other hand $43t^{3}-39t+32 \geq 40t^{3}-40t+30 = 10(4t^{3}-4t+3)$. The map $t \mapsto 4t^{3}-4t+3$ is positive at $t=0$ and $t=1$. It has a critical point at $t=1/\sqrt{3}$ where its value is $3-\frac{8 \sqrt{3}}{9}>0$. Thus $f''\geq 0$ on $[0,1]$. 
\end{proof}

\subsection{Proof of $C(4,6)=1/3$}
Let $p=4$ and $q=6$. We have 
\begin{align*}
 &\frac{(c^{q-2}(a^{q}-1)+c^{q}(a^{2}-a^{q})+1-a^{2})^{\frac{p-2}{q-p}}}{(c^{p-2}(a^{p}-1)+c^{p}(a^{2}-a^{p})+1-a^{2})^{\frac{q-2}{q-p}}}\cdot  (1-c)(1-a)(1-ac) -\frac{1}{3} \\
 &=\frac{c(2c-1)a^{2}-(c+1)^{2}a+3c^{2}-c+2}{3(1+c)(1+a)(1+ac)}.
\end{align*}
It suffices to show $\inf_{a,c \in (0,1)} c(2c-1)a^{2}-(c+1)^{2}a+3c^{2}-c+2 =0$. Indeed, the map $r(a) = c(2c-1)a^{2}-(c+1)^{2}a+3c^{2}-c+2$ is decreasing in $a$ because the linear function  $r'(a) = 2c(2c-1)a-(c+1)^{2}$ satisfies $r'(0)<0$ and $r'(1)=2c(2c-1)-(c+1)^{2} = 3c^{2}-4c-1 = -3c(1-c)-c-1<0$. On the other hand $r(1)=(2c-1)^{2}\geq 0$ and it becomes equality at $c=1/2$.

\subsection{Proof of Corrollary~\ref{texnika1}}
Let $X = | \sum_{j=1}^{n} a_{j} \xi_{j}|$, where $\xi_{i}$, $i=1, \ldots, n$,  are i.i.d. random variables defined as in (\ref{biased1}), and $\sum_{j=1}^{n} a_{j}^{2}=1$.   We have $\| X\|_{2}=1$, and 
\begin{align*}
\mathbb{E} |X|^{4} = 3+(\mathbb{E}\xi_{1}^{4}-3)\sum_{j=1}^{n}a_{j}^{4} \geq \min\{3, \mathbb{E} \xi_{1}^{4}\} = 1+\min\left\{2, \frac{(2p-1)^{2}}{p(1-p)}\right\}.
\end{align*}
On the other hand if $\xi'_{j}$ is independent copy of $\xi_{j}$, $j=1,\ldots, n$, and $\varepsilon_{j}$, $j=1, \ldots, n,$ are i.i.d. symmetric $\pm 1$ Rademacher random variables  we have b a symmetrization argument 
\begin{align*}
\mathbb{E} |X|^{6} &= \mathbb{E} |\sum_{j=1}^{n}a_{j} \xi_{j}|^{6} \leq \mathbb{E} |\sum_{j=1}^{n}a_{j} (\xi_{j}-\xi'_{j})|^{6} =\mathbb{E}_{\xi, \xi'} \mathbb{E}_{\varepsilon} |\sum_{j=1}^{n}a_{j} \varepsilon_{j}(\xi_{j}-\xi'_{j})|^{6}\\
& \leq 15\mathbb{E}_{\xi, \xi'}  \left(\sum_{j=1}^{n}a^{2}_{j} (\xi_{j}-\xi'_{j})^{2}\right)^{3}\leq \frac{15}{p^{3}(1-p)^{3}}, 
\end{align*}
where we used Khinchin's inequality $\mathbb{E} |\sum \varepsilon_{i}b_{i}|^{6} \leq 15 \left(\sum b_{i}^{2}\right)^{3}$ with the sharp constant $15$. 
Thus Theorem~\ref{main2} applied with constant $C(4,6) = 1/3$ gives 
\begin{align*}
\|X\|_{1} \leq 1-\frac{\min\left\{4p^{3}(1-p)^{3}, (2p-1)^{4}p(1-p)\right\}}{45-3(p(1-p))^{3}}.
\end{align*}

\subsection{Proof of Corrollary~\ref{texnika2}}\label{datvla}

We need the following identity obtained in \cite{pi2}.  For any $f:\{-1,1\}^n \to \mathbb{R}$, we have 
\begin{equation}
 f(x) - \mathbb{E} f(x) = \int_{0}^{\infty}\frac{1}{\sqrt{e^{2t} - 1}}\mathbb{E}_\zeta \sum_{j=1}^n \delta_j(t) D_j f(x\zeta(t)) dt, \label{heat}
\end{equation}
where $x\zeta(t) = (x_1\zeta_1(t),...,x_n\zeta_n(t))$, and the $\zeta_i(t)$ are i.i.d. random variables with
\[\mathbb{P}\{\zeta_i(t) = \pm 1\} = \frac{1\pm e^{-t}}{2},\]
and   $\delta_i = \frac{\zeta_i(t) - \mathbb{E} \zeta_{i}(t)}{\sqrt{\mathrm{Var}(\zeta_{i}(t))}}$.  We have 
\begin{align*}
&\mathbb{E} |f(x)-\mathbb{E} f(x)| \leq \int_{0}^{\infty} \mathbb{E} \left|\sum_{j=1}^{n} \delta_{j}(t) D_{j}f(x \zeta(t)) \right| \frac{dt}{\sqrt{e^{2t}-1}}\\
&=\int_{0}^{\infty}\mathbb{E} \left|\sum_{j=1}^{n} \delta_{j}(t) D_{j}f(x) \right| \frac{dt}{\sqrt{e^{2t}-1}}  \\
&\stackrel{(p=\frac{1+e^{-t}}{2})}{=} \int_{1/2}^{1}\mathbb{E} \left|\sum_{j=1}^{n} \xi_{j} D_{j}f(x) \right| \frac{dp}{\sqrt{p(1-p)}}  \stackrel{(\ref{dax1})}{\leq} \\
& \mathbb{E} |\nabla f| \int_{1/2}^{1}\left( 1-\frac{\min\left\{4p^{3}(1-p)^{3}, (2p-1)^{4}p(1-p)\right\}}{45-3(p(1-p))^{3}}\right) \frac{dp}{\sqrt{p(1-p)}}\\
&\approx \left(\frac{\pi}{2}-\delta\right) \mathbb{E} |\nabla f|,
\end{align*}
where 
$$
\delta = \int_{1/2}^{1}\frac{\min\left\{4p^{3}(1-p)^{3}, (2p-1)^{4}p(1-p)\right\}}{45-3(p(1-p))^{3}} \frac{dp}{\sqrt{p(1-p)}} \approx 0.00013...
$$

\subsection*{Acknowledgments} We thank Ramon van Handel for helpful comments. P.I was supported in part by NSF CAREER-DMS-2152401.


\end{document}